\documentclass[12pt,a4paper]{amsart}
\usepackage{amsfonts}
\usepackage{amsthm}
\usepackage{amsmath}
\usepackage{amscd}
\usepackage[latin2]{inputenc}
\usepackage{t1enc}
\usepackage[mathscr]{eucal}
\usepackage{indentfirst}
\usepackage{graphicx}
\usepackage{graphics}
\usepackage{pict2e}
\usepackage{epic}
\numberwithin{equation}{section}
\usepackage[margin=2.9cm]{geometry}
\usepackage{epstopdf} 
\usepackage{tikz}
\usepackage{tikzpagenodes}
\usepackage{gensymb}
\usepackage{amsmath}
\usepackage{stmaryrd}

\theoremstyle{plain}

 \theoremstyle{definition}

\newtheorem{?}[Th]{Problem}

\begin{document}

\title{A\MakeLowercase{n intuitive proof of the} D\MakeLowercase{voretzky-}H\MakeLowercase{anani theorem in} $\mathbb{R}^2$}

\author{
  E\MakeLowercase{fstratios}, M\MakeLowercase{arkou}\\
  U\MakeLowercase{niversity of} C\MakeLowercase{ambridge}\\
  \texttt{\MakeLowercase{em626@cam.ac.uk}}
}




\begin{abstract} 

The Dvoretzky-Hanani theorem states that the general term of any perfectly divergent series in a finite dimensional space does not tend to zero. An intuitive proof is provided $\mathbb{R}^2$,using a construction that allows us to determine a choice of $\pm$ such that

\begin{gather*}a_1 \pm a_2 \pm a_3 \pm a_4...\pm a_n ... \end{gather*}

\noindent{converges to a point in the space if $||a_i|| \to 0$. Extensions to the construction are proposed for the general $\mathbb{R}^n$.}

\end{abstract}

\maketitle
\pagestyle{plain}

\vspace{-0.5cm}
\section{Sketch of proof}
\vspace{0.75cm}
\noindent{For a sequence of vectors $\left\{ a_i \right\} \in \mathbb{R}^2$ let us call a particular choice of $\pm$ in}

\[a_1 \pm a_2 \pm a_3 \pm a_4 ...\pm a_n ...\]

\noindent{an \textit{assignment}, denoted $A$. For example}

\[A_1(\left\{ a_i \right\}) = a_1 + a_2 + a_3 + a_4 ...  + a_n ...\]
\[A_2(\left\{ a_i \right\}) = a_1 - a_2 + a_3 - a_4 ...  + (-1)^n a_n ...\]

\noindent{are two assignments for $\left\{ a_i \right\}$. We will partition the sequence $\left\{ a_i \right\} \in \mathbb{R}^2$ ($||a_i|| \to 0$) into finite subsets $S_m$ such that they preserve the order of the sequence $\left\{ a_i \right\}$, and also $a_i \in S_m \implies ||a_i|| <1/(m+1)^2$. We then provide $A_m$ such that}
\\
\begin{gather*}||A_m(\left\{ S_m \right\})|| < 6/(m+1)^2\end{gather*}
\\
As the terms $||A_m(\left\{ S_m \right\})||$ have size of at most $6/(m+1)^2$, using the $A_m$'s in summing $\left\{ a_i \right\}$ gives a sum of vectors which converges, proving the result.\\

\pagebreak

\section{Two-dimensional case}
\vspace{0.75cm}
\noindent{Let $\left\{ a_i \right\} \in \mathbb{R}^2$ with $||a_i|| \to 0$, under the Euclidian norm. As $||a_i|| \to 0$, there exists a minimum positive integer  $N_k \in \mathbb{Z^+}$ for any $k \geq 0 \in \mathbb{Z^+}$ such that}\\

\[||a_n|| <1/(k+1)^2, \forall n > N_k\]\\

Let $\{N_i\}$ be the sequence of all such indices and define $S_m = \left\{a_{N_m+1}, a_{N_m+2}, ..., a_{N_{m+1}} \right\}$ for $m > 1$ and $S_0 = \left\{a_{1}, a_{2}, ..., a_{N_{0}} \right\}$ for $m = 0$ respectively. Because $S_0$ is finite, its assignment  $A_0$ does not affect the convergence of the overall series. Under these definitions if $s \in S_m$, $m > 0$ then

\[1/(m+2)^2<||s||<1/(m+1)^2\]\\

\begin{figure}[b]
\begin{tikzpicture}[scale=0.6, every node/.style={transform shape}]

\draw (-5,0) -- (5,0);
\draw (0,-5) -- (0,5);
\draw (-0.5*5,-1.73205081/2*5) -- (0.5*5,1.73205081/2*5);
\draw (-0.5*5,1.73205081/2*5) -- (0.5*5,-1.73205081/2*5);
\draw (0,-5) -- (0,5);
\draw (1,0) arc (0:60:1);
\draw (0.5/3*2,1.73205081/3) arc (60:120:2/3);
\draw (-1,0) arc (180:120:1);
\draw (2/3,0) arc (0:-60:2/3);
\draw (-0.5,-1.73205081/2) arc (240:300:1);
\draw (-2/3,0) arc (180:240:2/3);
\node [label=below:$60^{\circ}$,draw,inner sep=0pt,minimum size=0pt] at (0.5*5/2,1.73205081/2*5/4) {};
\node [label=below:$60^{\circ}$,draw,inner sep=0pt,minimum size=0pt] at (0 , 1.25) {};
\node [label=below:$60^{\circ}$,draw,inner sep=0pt,minimum size=0pt] at (-0.5*5/2,1.73205081/2*5/4) {};
\node [label=below:$60^{\circ}$,draw,inner sep=0pt,minimum size=0pt] at (-0.5*5/2,-1.73205081/2*5/4+0.5) {};
\node [label=below:$60^{\circ}$,draw,inner sep=0pt,minimum size=0pt] at (0 , -1.1) {};
\node [label=below:$60^{\circ}$,draw,inner sep=0pt,minimum size=0pt] at (0.5*5/2,-1.73205081/2*5/4+0.5) {};
\node [label=below:$U_1$,draw,inner sep=0pt,minimum size=0pt] at (0.5*5,1.73205081/2*5/2) {};
\node [label=below:$U_2$,draw,inner sep=0pt,minimum size=0pt] at (0 , 1.25*5/2) {};
\node [label=below:$U_3$,draw,inner sep=0pt,minimum size=0pt] at (-0.5*5,1.73205081/2*5/2) {};
\node [label=below:$U_4$,draw,inner sep=0pt,minimum size=0pt] at (-0.5*5,-1.73205081/2*5/2+0.5) {};
\node [label=below:$U_5$,draw,inner sep=0pt,minimum size=0pt] at (0 , -1.25*5/2+0.1) {};
\node [label=below:$U_6$,draw,inner sep=0pt,minimum size=0pt] at (0.5*5,-1.73205081/2*5/2+0.5) {};

\end{tikzpicture}
\caption{Partitioning of $\mathbb{R}^2$ in six areas.}
\label{fig:areas}
\end{figure}
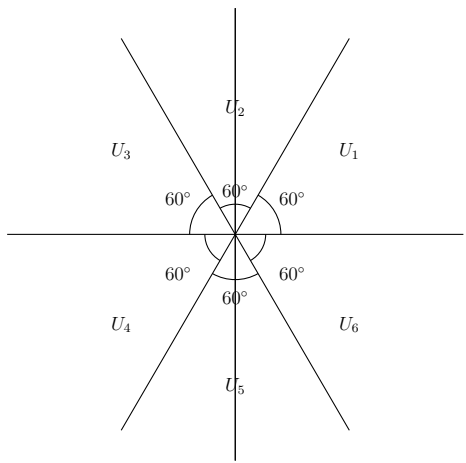

To provide $A_m$, we first partition $\mathbb{R}^2$ into six equal regions $U_i, 1\leq i \leq 6$ (see Fig. 1). If $S_m$ contains up to $6$ terms, any choice of $A_m$ will satisfy\\

\[||A_m(\left\{ S_m \right\})|| < 6/(m+1)^2\]\\

so we may consider  $S_m$ to have more than 6 terms. In this case we form distinct pairs of terms $(u, v)$ $u, v \in S_m, U_i$ as follows: start pairing the terms randomly and if any one remains unpaired leave it unpaired. Then calculate the difference $u - v$ for all pairs $(u, v)$. By our partitioning of $\mathbb{R}^2$, the angle between $u$ and $v$ is $\leq 60^{\circ}$ so by the cosine rule either $||u-v|| \leq ||u||$ or $||u-v|| \leq ||v||$ so

\begin{gather*}||u-v||<1/(m+1)^2\end{gather*}\\

The set $C_{m1}$ containing the $(u - v)$ terms and any unpaired terms satisfies $n(C_{m1}) < n(S_k)$ and $c \in C_{m1} \implies ||c|| < 1/(m+1)^2$. Repeating this process using $C_{m1}$ instead of $S_m$ gives another set $C_{m2}$ with $n(C_{m2}) < n(C_{m1})$ and $c \in C_{m2} \implies ||c|| < 1/(m+1)^2$. We can repeat the process for $i = \left\{3, 4, 5...\right\}$ to obtain gives $C_{mi}$'s with a strictly decreasing number of elements until for some $j$ $n(C_{mj}) \leq 6$ for some $j$.\\

Adding the (at most six) terms of $C_{mj}$ gives a vector with magnitude $\leq6/(m+1)^2$, and we may trace the above process backwards to determine $A_m$ such that\\

\[||A_m(\left\{ S_m \right\})|| < 6/(m+1)^2\]\\

and using these $A_m$'s in summing $\left\{ a_i \right\}$\\

\[A_0(\left\{ S_0 \right\}) + A_1(\left\{ S_1 \right\}) + A_2(\left\{ S_2 \right\}) + ... + A_m(\left\{ S_m \right\}) + ... \]\\

gives a sum of vectors which converges to a point, because each coordinate of 

Equivalently, because $||a_i|| \to 0$ implies the existence of a $A$ such that $A(\{a_i\})$ converges, which we can determine through this algorithm, a perfectly divergent series must require that $||a_i||$ does not tend to 0.\\

\section{Argument extensions}

By the above argument, for a non-Euclidian norm $||\cdot||_d$ over $\mathbb{R}^2$ with the property $||u+\delta u||_d \to ||u||_d$ as $||\delta u||_d \to 0, \forall u, \delta u \in \mathbb{R}^n$, there clearly also exists an assignment such that $A(\{a_i\})$ converges: algorithm for a Euclidian norm we can determine A such that $A(\left\{ a_i \right\}) \to l \in \mathbb{R}^2$ and $||A(\left\{ a_i \right\})||_d$ converges regardless of $||\cdot||_d$.\\

In order to extend the result to $\mathbb{R}^n$ with a Euclidian norm it suffices to prove that there exists a partitioning of $\mathbb{R}^n$ into a finite number of subsets such that for any $u, v \in \mathbb{R}^n$, either $||u-v|| \leq ||u||$ or $||u-v|| \leq ||v||$. The rest of the argument remains unchanged.

\end{document}